%% file: floating_sculpture_arxiv.tex
\documentclass[reqno]{amsart}
\usepackage{amsmath,amssymb,amsthm,graphicx,a4wide,enumerate}
\usepackage{mathtools,color,hyperref}
\usepackage[small,bf]{caption}
\theoremstyle{plain}

\theoremstyle{definition}

\theoremstyle{remark}

\newcommand{\prn}[1]{\left(#1\right)}
\newcommand{\abs}[1]{\left|#1\right|}
\newcommand{\brk}[1]{\left[#1\right]}

\newcommand{\norm}[1]{\left\|#1\right\|}

\newcommand{\ud}[1]{\,\text{d}#1}
\newcommand{\LBO}[1]{\mathcal{L}_{_{#1}}}

\begin{document}
\parskip1ex

\title[The Sound of an Evolving Floating Sculpture]
{The Sound of an Evolving Floating Sculpture}
\author{Benjamin Seibold}
\address[Benjamin Seibold]
{Department of Mathematics \\ Temple University \\
1805 North Broad Street \\ \newline Philadelphia, PA 19122}
\email{seibold@temple.edu}
\urladdr{http://www.math.temple.edu/\~{}seibold}
\author{Yossi Farjoun}
\address[Yossi Farjoun]
{G.\ Mill\'an Institute of Fluid Dynamics \\ Nanoscience and Industrial Mathematics \\
University Carlos III de Madrid \\ Av.\ Universidad 30, 28911 Legan\'es \\ Spain}
\email{yfarjoun@ing.uc3m.es}
\date{January 15, 2012}
\begin{abstract}
Commissioned by MIT's in-house artist Jane Philbrick, we evolve an abstract 2D surface (resembling Marta Pan's 1961 ``Sculpture Flottante I'') under mean curvature, all the while calculating the eigenmodes and eigenvalues of the Laplace-Beltrami operator on the resulting shapes. These are then synthesized into a sound-wave embodying the ``swan song'' of the surfaces as the evolve to points and vanish. The surface is approximated by a triangulation, and we present a robust approach to approximate the normal directions and the mean curvature. The resulting video and sound-track were parts in the Jane Philbrick's exhibition ``Everything Trembles'' in Lund, Sweden, 2009.
\end{abstract}

\maketitle

\section{Introduction}
This paper presents a project in the areas of Computation and Applied Mathematics. With both authors being researchers in Applied Mathematics and Scientific Computing, this does not sound unusual. However, the area of application \emph{is} quite unusual (at least for the authors), namely: visual art. Employing and devising a variety of methods in computational geometry, numerical linear algebra, and spectral methods, two pieces of art---an animation and a sounds piece---were created. They were part of the exhibition ``Everything Trembles'' that ran at the Skissernas Museum (Lund, Sweden) in fall 2009. In this paper, we present the mathematical modeling that was required to formulate the artistic problem, the computational approaches undertaken, and the steps from computational results to creating pieces of art.

In March 2009, through our colleague Enno Lenzmann (now University of Copenhagen), the visual artist Jane Philbrick (MIT Center for Advanced Visual Studies) approached the authors, with the following idea. She was preparing an exhibition fully devoted to Marta Pan's ``Sculpture Flottante I'', which is an exhibit in the Kr{\"o}ller-M{\"u}ller Museum (Otterlo, Netherlands). Her goal was to bring the sculpture to ``life'' in drawings, painting, and sculptures, as well as in modern visual art techniques, such as photographs, video, and sound installations. The specific question that she addressed us with was whether one can predict the sculpture's ``mathematical future'', by evolving and deforming it, as well as its characteristic sound.

Not unexpectedly, our first step was to ``model'' the artist's idea into more precise mathematical concepts. Our answer, then discussed and agreed on by the artist, was to let the sculpture (as defined by its 2D surface) evolve under a mean curvature flow (thus answering the question about the future of the sculpture), and to create sound-form based on eigenfunctions of the Laplace-Beltrami operator on the evolving surface. Both ideas move the sculpture away from its rigid shape, and instead interpret it as a flexible membrane with tension, very much like a punctured soap-bubble-like surface. The mean curvature flow pulls the surface together, creates pinch-off events that break the object apart, and eventually leads to its complete vanishing. The Laplace-Beltrami operator acts analogous to the classical Laplacian, but on a curved surface. Its eigenfunctions are used to evolve wave forms in time that are in turn used to create a characteristic sound pattern that changes as the surface evolves under the mean curvature flow. The two pieces of art created out of these mathematical concepts are: a video showing the evolution of the surface under mean curvature flow, and a sound file created out of the characteristic sounds of the surface as it evolves under the flow.

This paper is organized as follows. In Sect.~\ref{sec:mathematical_background}, an outline of the mathematical background relevant to this project is given. The ``transformation'' of the original sculpture into data that can be used for computational approaches is described in Sect.~\ref{sec:sculpture_to_computer_model}. The methods used and devised to compute the mean curvature flow are presented in Sect.~\ref{sec:geometry_evolution}, and the approaches to calculate the characteristic sounds of the surface are given in Sect.~\ref{sec:surface_sounds}. Finally, in Sect.~\ref{sec:results_to_art}, we describe how the computational results are assembled into pieces of art.

\section{Mathematical Background}
\label{sec:mathematical_background}

\subsection{Curvature}
In the following, we consider smooth, simple curves in $\mathbb{R}^2$ and smooth, orientable, and closed surfaces in $\mathbb{R}^3$. At any point on the curve/surface, $\hat{n}$ is the outward normal vector. For a smooth curve, $C$, embedded in $\mathbb{R}^2$, the curvature $\kappa$ at a point $P$ has an absolute value equal to the reciprocal of the radius of the osculating circle, i.e.~the circle that passes through $P$ and a pair of additional points on the curve that lie infinitesimally close to $P$ \cite{Yates1952}. Moreover, we use the convention that $\kappa$ is negative if $\hat{n}$ points towards the center of the osculating circle, and positive otherwise. With this convention, a circle of radius $r$ has a constant curvature $\kappa = \frac{1}{r}$.

For a smooth surface, $S$, embedded in $\mathbb{R}^3$, multiple notions of curvature exist. At any point, $P$, on $S$, one can consider all curves on $S$ that pass through $P$, and find the corresponding curvature for each of them (using the concepts established above). Of those curvatures, the maximal, $\kappa_1$, and the minimal, $\kappa_2$, are called the ``principal curvatures'' of $S$ at $P$. Their product, $\kappa_1\kappa_2$, is called ``Gaussian curvature'', and their average, $\kappa = \frac{1}{2}(\kappa_1+\kappa_2)$, is called ``mean curvature'' \cite{Yates1952}.

While Gaussian curvature as an intrinsic quantity is independent of the surface's (isometric) embedding, mean curvature is extrinsic, i.e.~it depends on the embedding of the surface in its ambient space $\mathbb{R}^3$. Mean curvature is closely related to the physical effect of surface tension, which acts normal to the interface between two immiscible fluids and is proportional in magnitude to its mean curvature. The mathematical relevance of the mean curvature was first formulated in the Young-Laplace equation \cite{Young1805, Laplace1806}, which describes the scenario in which surface tension is balanced by a capillary pressure difference sustained across the interface between two static fluids.

When there is no pressure counter-acting the surface tension on an interface between fluids (or a thin film) will move it as to decrease to surface area. A mathematical model for such movements due to surface tension is mean curvature flow \cite{Ecker2004}, which is a geometric nonlinear partial differential equation \cite{Ecker2004}. According to this flow, a surface embedded in $\mathbb{R}^3$, or a curve embedded in $\mathbb{R}^2$, moves parallel to its normal vector with a velocity proportional to its curvature. Specifically, a point $\vec{x}$ on the surface/curve moves according to the law of motion
\begin{equation}
\label{eq:mean_curvature_flow}
\dot{\vec{x}} = -\kappa(\vec{x})\hat{n}(\vec{x})\;,
\end{equation}
where $\kappa$ is the (mean) curvature, and $\hat{n}$ is the outward pointing normal vector, both evaluated at $\vec{x}$. The negative sign is due to our convention that curvature is positive where the surface is convex, and negative where it is concave.

Curvature flow of curves in $\mathbb{R}^2$ exhibits a relatively regular behavior. Any simple, closed curve shrinks continuously, and at a steady rate, until it eventually vanishes into a single point \cite{Grayson1987}. Due to the Gauss-Bonnet theorem \cite{GuilleminPollack1974}, the total curvature over the curve is an invariant, i.e.~$\int_{C}\kappa\ud{S} = 2\pi$, and thus the rate of change of area enclosed by the curve equals
$\frac{d}{dt}A(t) = \int_{C}\dot{\vec{x}}\cdot\hat{n}\ud{S} = -\int_{C}\kappa\ud{S} = -2\pi$. Hence, the time of the curve's vanishing is $T = \frac{A(0)}{2\pi}$, where $A(0)$ is the area enclosed by the initial curve. To the authors' knowledge, it is yet an unsolved problem to express the exact point of the curve's vanishing in terms of geometric properties of its initial configuration \cite{Sethian1989}.

Mean curvature flow of surfaces in $\mathbb{R}^3$ is far less regular than curvature flow of curves. First, the Gauss-Bonnet theorem does not apply to mean curvature (it applies to Gaussian curvature, which is not considered here). Second, a surface that moves under mean curvature flow can separate into multiple disjoint components. Separation happens in a ``pinch-off'' event, in which the neck of a dumbbell-shaped part of the surface shrinks to radius zero. Since mean curvature flow is closely related to surface tension of physical fluids, this pinch-off phenomenon is well-known in fluid dynamics as a Plateau-Rayleigh instability, in which a falling stream of fluid breaks up into smaller droplets \cite{Plateau1873, Rayleigh1879}. In fact, in the evolution of  ``Sculpture Flottante I'', we observe such pinch-off events.

\subsection{Laplace-Beltrami Operator}
\label{subsec:laplace_beltrami}
Since we are also interested in the ``characteristic sound of the sculpture'', we need to consider modes of vibration of a surface in $\mathbb{R}^3$. In Sect.~\ref{sec:surface_sounds}, we use the surface wave equation
\begin{equation}
\label{eq:wave_equation}
u_{tt} = \LBO{} u
\end{equation}
as a model for sound wave propagation on the surface, where $\LBO{}$ is the Laplace-Beltrami operator, which is the generalization of the standard Laplace operator to curved surfaces. It is defined as
\begin{equation*}
\LBO{} u = \frac{1}{\sqrt{|g|}} \sum_i \partial_i\prn{\sqrt{|g|}
\sum_j g^{ij} \partial_j u}\;,
\end{equation*}
where $g$ is the metric tensor associated with the surface, and $|g| = |\det(g_{ij})|$ is the absolute value of its determinant. Since the Laplace-Beltrami operator encodes information about the surface's metric tensor, there is a close relation between it and mean curvature, namely (with our sign convention for curvature)
\begin{equation}
\label{eq:relation_LB_curvature}
\LBO{}\vec{x} = -2\vec{\kappa}(\vec{x})\;,
\end{equation}
where $\vec{x}$ is the position vector of a point on the surface, and $\vec{\kappa} = \kappa\hat{n}$ is the mean curvature normal vector \cite{Willmore1993} at the same position. As a consequence, the equation of mean curvature flow \eqref{eq:mean_curvature_flow} can also be written in terms of the Laplace-Beltrami operator as
\begin{equation*}
\dot{\vec{x}} = \tfrac{1}{2}\LBO{}\vec{x}\;,
\end{equation*}
which illustrates the similarity of mean curvature flow to a diffusion process: ripples in the surface shape are smeared out at a rate that is faster for ripples with higher spacial frequency.

The surface wave equation \eqref{eq:wave_equation} is a linear partial differential equation, whose solutions can be conveniently expressed in terms of eigenfunctions of the Laplace-Beltrami operator. Clearly, these eigenfunctions relate to important questions about the vibration of objects, such as Mark Kac's fundamental question ``Can one hear the shape of a drum?''\footnote{No.} \cite{Kac1966}. Moreover, the spectrum of the Laplace-Beltrami operator plays an essential role in many applications that go beyond the study of sound waves. Examples include signal processing and mesh design \cite{Taubin1995}, shape matching \cite{ReuterWolterPeinecke2005}, geometry processing \cite{Levy2006}, shape analysis and image segmentation \cite{ReuterBiasottiGiorgiPataneSpagnuolo2009}.

\subsection{Computational Geometry}
\label{subsec:computational_geometry}
In order to implement the mathematical concepts described above, specifically mean curvature flow and eigenmodes of the Laplace-Beltrami operator, on a computer, these geometric operators must be approximated on a surface that can be described by a finite amount of data. We consider the surface be given by a regular triangulation. Hence, we must employ approaches to approximate surface normals, curvature, and the Laplace-Beltrami operator on such a discrete geometry. The first approaches of approximating differential operators on triangulations go back to Delaunay \cite{Delaunay1934}. Since then, a wide variety of methods has been presented. However, to the authors' knowledge, the quote ``Despite extensive use of triangle meshes in Computer Graphics, there is no consensus on the most appropriate way to estimate simple geometric attributes such as normal vectors and curvatures on discrete surfaces'' \cite{MeyerDesbrunSchroderBarr2002} still holds true nowadays. One key challenge is to devise efficient approximation approaches that work in a robust fashion under general circumstances. Below, we outline a few representative approaches for the approximations of curvature and the Laplace-Beltrami operator.

Many approaches to approximate the Laplace-Beltrami operator at a vertex of a triangulation are based on so-called ``cotan formulas''. Let $V$ denote the set of vertices, and $E$ the set of pairs of vertices that are connected by an edge. Assume that a function $u$ is given on the vertices. Then $\LBO{}u$ at a vertex $x_i$ is approximated by
\begin{equation}
\label{eq:discrete_Laplace_Beltrami}
\LBO{} u(\vec{x}_i) \approx \hspace{-1.5em}
\sum_{\vec{x}_j\in V:(\vec{x}_i,\vec{x}_j)\in E} \hspace{-1.5em}
w_{ij} (u(\vec{x}_j)-u(\vec{x}_i))\;,
\end{equation}
where $w_{ij}$ are appropriate weights. A common choice for the weights \cite{BobenkoSpringborn2007} is given by
\begin{equation}
\label{eq:LB_weights}
w_{ij} = \tfrac{1}{2}(\cot\alpha_{ij}+\cot\alpha_{ji})\;,
\end{equation}
where $\alpha_{ij}$ and $\alpha_{ji}$ are the angles that are opposite to the edge $(\vec{x}_i,\vec{x}_j)$, in the two triangles adjacent to this edge (see, for example, Fig.~1 in \cite{BobenkoSpringborn2007}). Other, similar formulas, albeit with different normalization constants, have been studied in the literature \cite{DesbrunMeyerSchroderBarr1999, MeyerDesbrunSchroderBarr2002}.

Due to relation \eqref{eq:relation_LB_curvature}, an approximation rule for the Laplace-Beltrami operator induces an approximation to the mean curvature normal vector $\vec{\kappa} = \kappa\hat{n}$, by applying the formula to the surface coordinate vector $\vec{x}$. For instance, formula \eqref{eq:discrete_Laplace_Beltrami} leads to
\begin{equation}
\label{eq:curvature_cotan}
\vec{\kappa} \approx \hspace{-1.5em}
\sum_{\vec{x}_j\in V:(\vec{x}_i,\vec{x}_j)\in E} \hspace{-1.5em}
w_{ij} (\vec{x}_i-\vec{x}_j)\;.
\end{equation}
In words: an approximate curvature normal vector at a vertex $x_i$ is obtained as a weighted sum of all adjacent edge vectors that point into $x_i$. In many cases, cotan formulas yield good quality approximations to mean curvature, but not always: for example, this method does not perform well on extremely pointy configurations \cite[Sect.~4.1]{DesbrunMeyerSchroderBarr1999}. Another shortcoming of cotan formulas is that they do not reveal the full curvature tensor. Examples of alternative approaches include the ``Functional Laplacian'' \cite{BelkinSunWang2008}, which is based on Green's functions, and approximations using the Willmore energy of surfaces \cite{HildebrandtPolthier2011}. An alternative way to approximate curvature is least-squares fitting for principal directions \cite{MeyerDesbrunSchroderBarr2002}.

In this project, we employ the classical formula \eqref{eq:discrete_Laplace_Beltrami} for the computation of the spectrum of Laplace-Beltrami operator (see Sect.~\ref{sec:surface_sounds}), and a non-classical way to approximate curvature (based on least-square minimization) to achieve a stable and robust approximation of the mean curvature flow. This latter approach is described in detail in Sect.~\ref{subsec:approximation_normals_curvature}.

It should be remarked that there are other computational methods for mean curvature flows that are not based on explicit representations (e.g.~triangulations) of the surface. An important example is the level set method \cite{OsherSethian1988}, which represents the surface implicitly as the zero contour of a level set function $\phi$. One then can formulate a partial differential equation (in the ambient space) for $\phi$ so that its zero contour moves according to mean curvature flow \cite{Sethian1989, ChoppSethian1993}. Advantages of the level set approach are the conceptual simplicity of solving a PDE in the ambient space, and the automatic treatment of topology changes. Disadvantages are the extra problem dimension (due to moving to the surface's ambient space) and the loss of sharp surface features due to numerical viscosity. In addition, if the original configuration is given as a triangulation (as it is in our case), level set approaches would require back-and-forth transformations between explicit and implicit representations, each of which would incur errors.

\section{Transformation of the Sculpture into a Computer Model}
\label{sec:sculpture_to_computer_model}
\begin{figure}
\begin{minipage}[b]{.4050\textwidth}
\includegraphics[width=\textwidth]{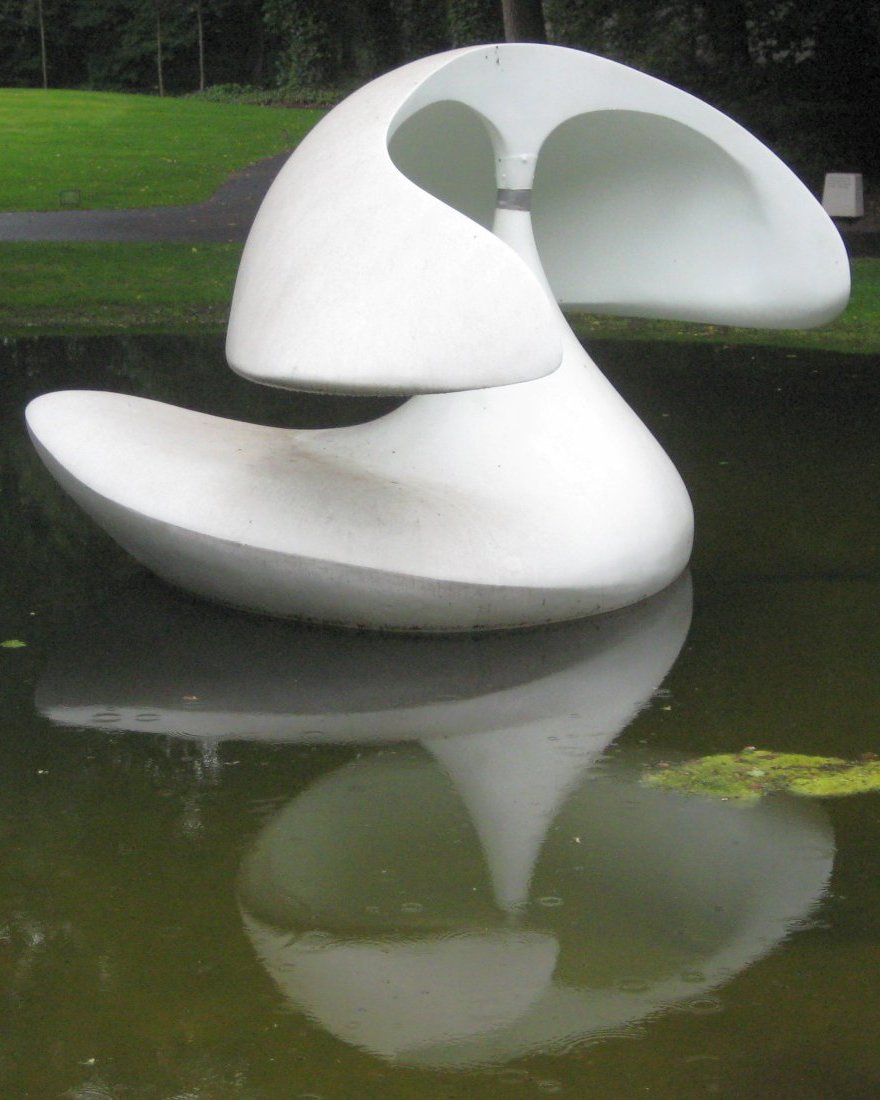}
\caption[Marta Pan's ``Sculpture Flottante I'' at the Kr{\"o}ller-M{\"u}ller Museum.]{Marta Pan's ``Sculpture Flottante I'' at the Kr{\"o}ller-M{\"u}ller Museum.\footnotemark\addtocounter{footnote}{-1}}
\label{fig:sculpture_flottante}
\end{minipage}
\hfill
\begin{minipage}[b]{.2235\textwidth}
\includegraphics[width=\textwidth]{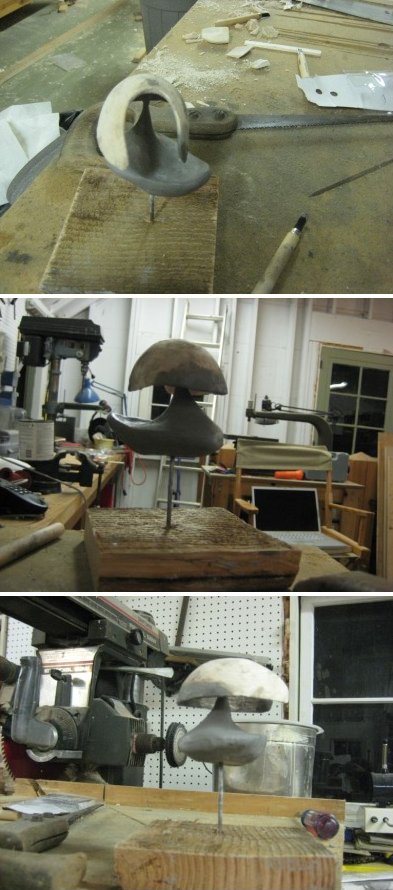}
\caption[Creation of a clay model.]{Creation of a clay model.\footnotemark\addtocounter{footnote}{-1}}
\label{fig:clay_model_creation}
\end{minipage}
\hfill
\begin{minipage}[b]{.3315\textwidth}
\includegraphics[width=\textwidth]{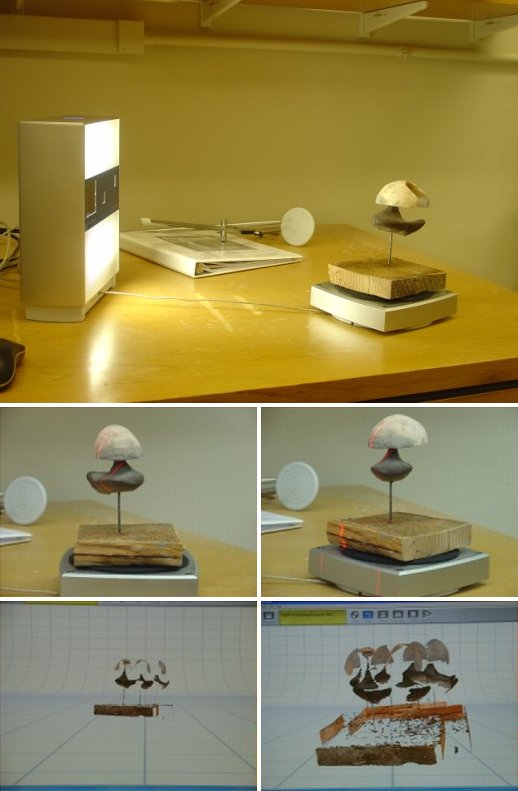}
\caption[Laser scanning of the model.]{Laser scanning of the model.\footnotemark}
\label{fig:clay_model_scan}
\end{minipage}
\end{figure}

The ``Sculpture Flottante I'' is a permanent exhibit in the sculpture garden of the Kr{\"o}ller-M{\"u}ller Museum (Otterlo, The Netherlands). It was created in 1961 by Marta Pan (1923--2008), a French sculptor of Hungarian birth. The sculpture, shown in its installation view in Fig.~\ref{fig:sculpture_flottante}, is made of fiberglass reinforced polyester. It consists of a base which is floating on a pond, and a top which is shaped like a cut-open umbrella and is connected to the base by a thin shaft.

In order to admit computational tasks with the sculpture's surface (as well as its visualization), a computer model needed to be generated. Since a direct scanning of the original sculpture turned out not possible, a clay replica of the sculpture was created (shown in Fig.~\ref{fig:clay_model_creation}), based on a series of photographs of the ``Sculpture Flottante I''. The clay model was then laser scanned (shown in Fig.~\ref{fig:clay_model_scan}) resulting in a stereolithography CAD file (\texttt{.stl} format), from which a triangulation of the 2D surface that resembled the sculpture could be extracted. This triangulation, which consists of more than 500,000 elements, was then transformed into a format that allowed us to perform computational geometry operations, as described in Sect.~\ref{subsec:computational_geometry}, Sect.~\ref{sec:geometry_evolution}, and Sect.~\ref{sec:surface_sounds}.

\footnotetext{All photographs courtesy of Jane Philbrick (MIT Center for Advanced Visual Studies).}

As can be seen in Fig.~\ref{fig:clay_model_creation}, the true sculpture's top and bottom are connected through a very thin shaft. Mean curvature flow contracts this shaft almost immediately, resulting is a pinch-off that separates the top from the bottom. We choose the time of this pinch-off as the starting time of the animation, described in Sect.~\ref{sec:results_to_art}.

\section{Evolution of the Geometry under Mean Curvature Flow}
\label{sec:geometry_evolution}
The computational task is to let the initial 2D surface evolve under mean curvature flow \eqref{eq:mean_curvature_flow}, until it has vanished completely. The basic idea is to advance equation \eqref{eq:mean_curvature_flow} forward in time using an explicit forward Euler scheme. This amounts to the following update scheme. In each step:
\begin{enumerate}[i.]
\item Given the triangulation, evaluate the mean curvature normal vector in the right hand of \eqref{eq:mean_curvature_flow}.
\item Move each vertex of the triangulation by $\Delta t$ times this right hand side vector, where $\Delta t$ is an appropriately chosen time step (see below).
\item If required, modify the triangulation (see below).
\end{enumerate}
In order to make this basic approach work, several technical challenges must be addressed.

First, the mean curvature normal vector on the vertices of a triangulation must be approximated in a robust fashion, such that motion under \eqref{eq:mean_curvature_flow} results in a stable scheme. It must be stressed that the approximation of the curvature normal vector for the sake of an evolution equation is significantly more demanding than a single evaluation of these differential quantities on given data. As an analogy, consider the approximation of the first derivative of a function using central differencing. When done once, this is a good idea. However, when used with forward Euler time stepping to approximate the linear advection equation, an unstable numerical scheme results. The methodology employed here is described in Sect.~\ref{subsec:approximation_normals_curvature}.

Second, as particles move, the topology of the surface may change (in a pinch-off event). The computational approach must be able to mimic such topology changes. This is described in Sect.~\ref{subsec:vertex_removal}.

Third, as the surface shrinks, particles tend to move closer towards their neighbors, and as they get closer the largest admissible time step shrinks. To guarantee that a sufficiently large time step can be taken, we merge neighboring vertices if they become too close, and then check for possible resulting topology changes. This aspects is described in Sect.~\ref{subsec:vertex_removal}.

Finally, the aspect of computational efficiency must be addressed. The maximum admissible time step is governed by the smallest features of the complete geometry. As a consequence, in most other parts of the surface, this time step is unnecessarily small, and thus a full update of the curvature normal vector everywhere would be prohibitively costly. Our remedy, described in Sect.~\ref{subsec:staleness}, is a mechanism by which the normal vector and curvature are only calculated where and when they truly need to be re-calculated. As long as they do not require an update, the last calculated values are used. Effectively, this results in a scheme with spatially and temporally adaptive local time steps. Since the computation of the differential quantities is a costly part of the computation, this added adaptivity reduces the computation time dramatically.

\subsection{Approximation of Normal Vectors and Mean Curvature}
\label{subsec:approximation_normals_curvature}
As described in Sect.~\ref{subsec:computational_geometry}, there is a wide variety of methods to approximate surface normal vectors and curvature, some of the most popular approaches being the ``cotan formulas'' \eqref{eq:curvature_cotan}. Approximations based on these formulas can be used to compute mean curvature flows, at least for short times, as it is sufficient for surface smoothing (see \cite{DesbrunMeyerSchroderBarr1999} for an example). For the task of evolving a surface all the way to its vanishing, and including topology changes, we have found that an alternative approach, described below, is more robust. This approach finds the normal vectors in a first step (see Sect.~\ref{subsubsec:approximation_normals}), and then finds the curvature tensor using a least-squares fitting approach (see Sect.~\ref{subsubsec:approximation_curvature}). While this two-step approach turns out to yield a robust mean curvature flow evolution, we are admittedly overdoing things a bit in having access to the normal vector and the full curvature tensor. As a consequence, the approximations presented here have the potential to be used for geometric flows that are more complicated than mean-curvature flow---a topic that we do not follow up with in this paper.

\subsubsection{Approximation of normal vectors}
\label{subsubsec:approximation_normals}
\begin{figure}
\begin{minipage}{.48\textwidth}
\vspace{1.24em}
\includegraphics[width=\textwidth]{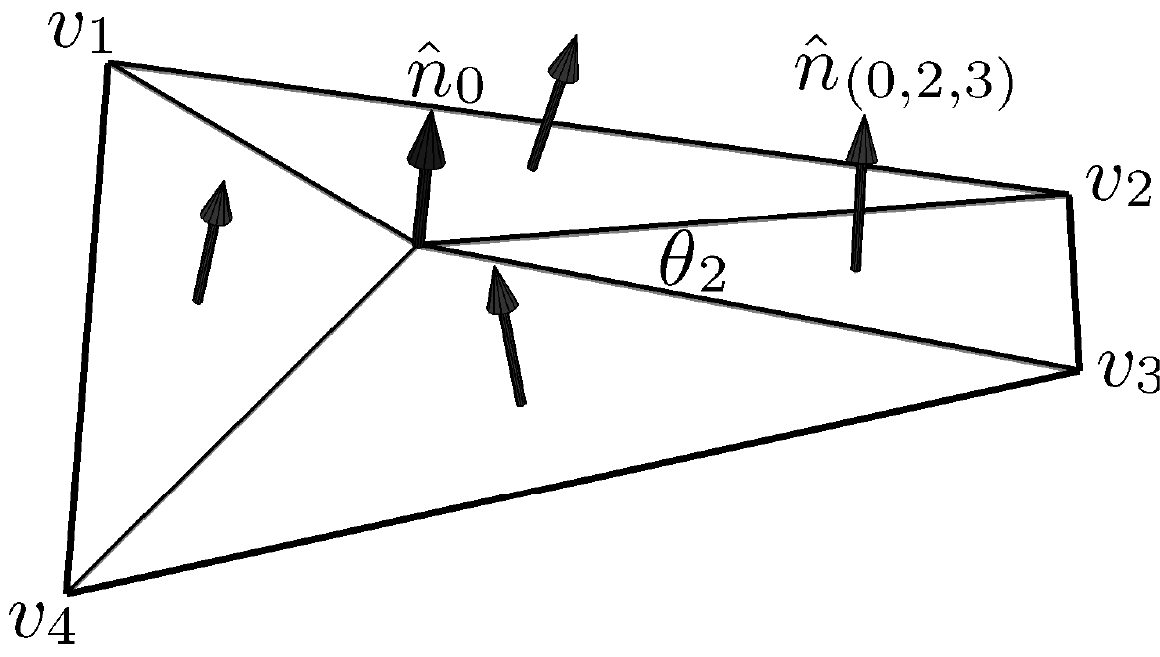}
\caption{The face-averaged normal vector.}
\label{fig:face:average}
\end{minipage}
\hfill
\begin{minipage}{.48\textwidth}
\includegraphics[width=\textwidth]{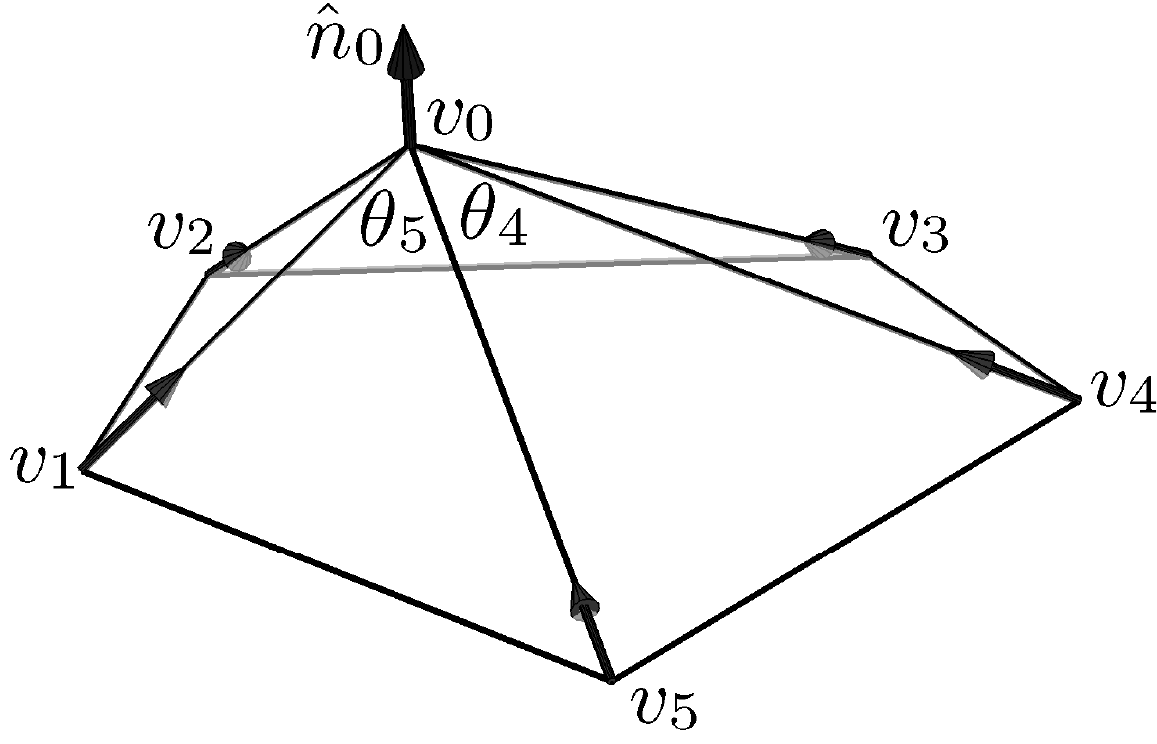}
\vspace{-1.5em}
\caption{The edge-averaged normal vector.}
\label{fig:edge:average}
\end{minipage}
\end{figure}

To approximate the normal vector at a vertex $v_0$, we look at two types of averages and chose the larger between them. In the following the center vertex $v_0$ is assumed to be at the origin as this simplifies the notation. Moreover, we consider the $i$ direct neighbor vertices $v_1,\dots,v_i$ be ordered counterclockwise around the center vertex when the shape is viewed from the outside. This is possible for all orientable surfaces. The choice of $v_0$ as the origin implies that the vertices $v_j$ are also the edges connecting the $v_0$ to the vertices.
\begin{itemize}
\item
The average face normals of adjacent faces around the vertex, weighted by the angle at the vertex (see Fig.~\ref{fig:face:average}):
\begin{align*}
\text{Face normal average}&=\frac{ \sum_{j=1}^i \theta_j \hat n_{(0,j,j+1)}}{\sum_{j=1}^i \theta_j}\;,
\shortintertext{where}
\theta_j&=\arcsin\prn{\frac{\norm{v_j\times
v_{j+1}}}{\norm{v_j}\norm{v_{j+1}}}}
\shortintertext{is the angle between the edges $j$ and $j+1$, and}
\hat n_{(0,j,j+1)}&=\frac{v_j\times v_{j+1}}{\norm{v_j\times v_{j+1}}}
\end{align*}
is the normal to the face that is between the $j$ and $j+1$ edges.
\item
The average of normalized adjacent \emph{edges}, pointing toward the center vertex, $-v_j$, weighted by the sum of the two angles at the vertex and that edge (see Fig.~\ref{fig:edge:average}):
\begin{align}
\label{eq:9}
\text{Edge direction average} &=-\frac{  \sum_{j=1}^i
\prn{\theta_{j-1} + \theta_{j} }\frac{v_j}{\norm{v_j}}}{2\sum_{j=1}^i \theta_j}\;.
\end{align}
Note that \eqref{eq:9} is similar in flavor to the cotan approximation \eqref{eq:curvature_cotan} to the curvature normal vector.
\end{itemize}
Of these two averages we use the one that is larger in magnitude, since this would imply that the parts being averaged agree on the normal to a greater extent. The reason we are using two different methods is that we need to find an accurate normal direction in both extremal cases: (a) when the faces are part of a large, nearly flat domain, and (b) when neighboring triangles are at the knife's edge of a surface, or are part of a very narrow spike. In the friendly case (a), the face normal average yields a good normal direction, since the face normals are almost parallel. In contrast, the edge direction average possesses only a small component in the direction normal to the surface, and is henceforth not used. In the case (b), when the center vertex is on a tip of a very sharp, needle-like, formation, the face normal average is very small (and thus discarded), while the edge direction average is nearly a unit vector in the same direction as the sharp formation.

When evolving a surface under mean curvature flow, in the vast majority of times the close-to-flat case (a) is observed, as the regularizing flow tends to smooth the surface. However, the proper treatment of the ``spiky'' case (b) is still required to obtain a robust treatment of the initial conditions, and of pinch-off events. Whichever type of average is used, the resulting vector is, of course, normalized before being used as the surface normal.

\subsubsection{Approximation of the curvature tensor}
\label{subsubsec:approximation_curvature}
Having found a normal vector $\hat n$ at the vertex $v_0$, we now need to approximate the curvature. To do so, we introduce a local coordinate system, where $\hat z = \hat n$ is that approximate normal, and $\hat x$ and $\hat y$ span the plane $\hat n^\perp$. The following formulae are written with regards to that coordinate system, where $v_0$ is the origin and the points $(x,\,y,\,z)$ refers to $v_0+ x\hat x+y\hat y +z\hat z$. In these coordinates the neighboring vertex $v_j$ is located at $(x_j,\,y_j,\,z_j)$.

We now approximate the curvature by fitting a quadratic polynomial
\begin{equation}
  \label{eq:quadratic:form}
a + \vec x^T
  \begin{pmatrix}
    d&e\\e&f
  \end{pmatrix}
 \vec x
\end{equation}
to the triangulated surface $\bigcup_{j=1}^i(v_0=\vec 0,v_j, v_{j+1}) $, where, for convenience of notation, we define $v_{i+1} = v_1$, and furthermore $\vec x = (x,y)^T$.
To approximate the shape constant (i.e.~the curvature tensor and the constant shift) of the surface, we look to minimize the functional
\begin{align}
  \label{eq:err:1}
E[\bold{w}]&=\sum_{j=1}^i \iint_{\text{triangle } j} \prn{ a + \begin{pmatrix}x\\y\end{pmatrix}^T
  \begin{pmatrix}
    d&e\\e&f
  \end{pmatrix}
 \begin{pmatrix}x\\y\end{pmatrix}
-z(x,y)}^2 \ud{x}\ud{y}\;,
\intertext{where $\bold{w}=\prn{a,\,d,\,e,\,f}^T$ is the vector of shape constants to the surface at a given vertex, and $z(x,y)$ is the $z$-value of the triangulation at the point $\prn{\begin{smallmatrix}x\\y\end{smallmatrix}}$. We re-write \eqref{eq:err:1} using a vector of functions,
$\bold{s}_{\prn{\begin{smallmatrix}x\\y\end{smallmatrix}}}=\prn{
  1,\,
    x^2,\,
   2xy ,\,
y^2}
^T$, as}
E[\bold{w}]&=\sum_{j=1}^i \iint_{\text{triangle } j} \prn{ \bold{s}^T_{\prn{\begin{smallmatrix}x\\y\end{smallmatrix}}}
\bold{w}
-z(x,y)}^2 \ud{x}\ud{y}\;. \nonumber
\intertext{In each triangle $j$ we use barycentric coordinates $\alpha$, $\beta\,$ and $\gamma$ (using the three vertices of the $j-$th triangle) where $\alpha+\beta+\gamma = 1$, to identify a point $\vec x$ inside the triangle:}
\begin{pmatrix}x\\y\end{pmatrix}
&=\begin{pmatrix}
      x_j& x_{j+1}\\
y_j&y_{j+1}
    \end{pmatrix}\begin{pmatrix}
     \alpha\\\beta
    \end{pmatrix}=X_j\begin{pmatrix}
     \alpha\\\beta
\end{pmatrix}\;.  \nonumber
\intertext{Due to the choice $v_0 = 0$, the value $\gamma$ does not appear in the expressions. Changing the integration variables to $\alpha$ and $\beta$, we obtain}
E[\bold{w}]&= \sum_{j=1}^i
\abs{\begin{matrix}
      x_j& x_{j+1}\\
y_j&y_{j+1}
    \end{matrix}}
\int_0^1\int_0^{1-\alpha}
\prn{ \bold{s}^T_{\prn{\begin{smallmatrix}
      x_j& x_{j+1}\\
y_j&y_{j+1}
    \end{smallmatrix}}\prn{\begin{smallmatrix}
     \alpha\\\beta
    \end{smallmatrix}}}
 \bold{w}
-\begin{pmatrix}
  \alpha\\\beta
\end{pmatrix}^T
\begin{pmatrix}
  z_j\\z_{j+1}
\end{pmatrix}}^2\,
\ud{\beta}\ud{\alpha}\;, \nonumber
\end{align}
where the determinant of the Jacobian enters due to the change of variables.

From the definition of $\bold{s}$ one can verify that applying a transformation
$X_j=\prn{\begin{smallmatrix} x_j& x_{j+1}\\
    y_j&y_{j+1} \end{smallmatrix}}$
to the argument of $\bold{s}$ is equivalent to right-multiplying by a matrix:
\begin{align}
\bold{s}_{X_j\vec x}=\bold{s}_{\vec{x}}
\begin{pmatrix}
1&0&0&0\\
0&x_j^2 & 2 x_j y_j & y_j^2 \\
0&  x_j x_{j+1} & x_{j+1} y_j+x_j y_{j+1} & y_j y_{j+1} \\
0&  x_{j+1}^2 & 2 x_{j+1} y_{j+1} & y_{j+1}^2
\end{pmatrix}=\bold{s}_{\vec{x}} \Xi_j\;.
\label{eq:second:symmetric:power}
\end{align}
The non-trivial lower $3\times 3$ block in the matrix in \eqref{eq:second:symmetric:power} is $S^2(X_j)$, the second symmetric power of $X_j$.
Letting $\vec{\alpha}=\prn{\begin{smallmatrix}
 \alpha\\\beta
 \end{smallmatrix}}$, and $\vec{z}_j=\prn{\begin{smallmatrix}
 z_j\\z_{j+1}
 \end{smallmatrix}}$,
the minimization functional can be written as:
\begin{equation*}
E[\bold{w}] = \sum_{j=1}^i
\abs{X_j}
\int_0^1\int_0^{1-\alpha}
\prn{ \bold{s}^T_{\vec{\alpha}}\Xi_j
\bold{w}-
\vec{\alpha}^T
\vec{z_j}}^2\,
\ud{\beta}\ud{\alpha}\;.
\end{equation*}
To integrate, we expand the power as a product of a vector with its transpose:
\begin{equation}
E[\bold{w}] = \sum_{j=1}^i
\abs{X_j}
\int_0^1\int_0^{1-\alpha}
\bold{w}^T
\Xi_j^T
\bold{s}_{\vec{\alpha}}
\bold{s}^T_{\vec{\alpha}}
\Xi_j
\bold{w}
-
2
\bold{w}^T
\Xi_j^T
\bold{s}_{\vec{\alpha}}
\vec{\alpha}^T
\vec{z_j}
+
\vec{z_j}^T
\vec{\alpha}
\vec{\alpha}^T
\vec{z_j}
\,
\ud{\beta}\ud{\alpha}.
\label{eq:last:integral}
\end{equation}
In \eqref{eq:last:integral} the only terms that depend on $\vec{\alpha}$ are $\bold{s}_{\vec{\alpha}} \bold{s}^T_{\vec{\alpha}}$, $\bold{s}_{\vec{\alpha}} \vec{\alpha}^T$, and $\vec{\alpha} \vec{\alpha}^T$. These terms' integrals are
\begin{equation*}
\iint \bold{s}_{\vec{\alpha}} \bold{s}^T_{\vec{\alpha}} = \tfrac{1}{180}
\begin{pmatrix}
 90 & 15 & 15 & 15 \\
 15 & 6 & 3 & 1 \\
 15 & 3 & 4 & 3 \\
 15 & 1 & 3 & 6
\end{pmatrix}
,\quad
\iint \bold{s}_{\vec{\alpha}} \vec{\alpha}^T = \tfrac{1}{60}
\begin{pmatrix}
 10 & 10 \\
 3 & 1 \\
 2 & 2 \\
 1 & 3
\end{pmatrix},\quad
\iint\vec{\alpha} \vec{\alpha}^T = \tfrac{1}{24}
\begin{pmatrix}
    2 & 1 \\
    1 & 2
\end{pmatrix},
\end{equation*}
which we label $A, B, C$ respectively.

To find where the minimum of $E[\bold{w}]$ is attained, we replace the integrals in \eqref{eq:last:integral} by $A$, $B$, and $C$ and differentiate with respect to $\bold{w}$.
This gives rise to the linear system for the minimizing $\bold{w}$:
\begin{equation}
  \prn{
  \sum_{j=1}^i
\abs{X_j}
 \Xi_j^T
 A
\Xi_j}
\bold{w}
=\sum_{j=1}^i
\abs{X_j}
\Xi_j^T
B \vec{z}_j\;.
\label{eq:linear:system:w}
\end{equation}
The mean curvature is now given by $\kappa = \frac{d+f}{2} = \prn{0,\frac{1}{2},0,\frac{1}{2}}\cdot\bold{w}$. Note that according to our approximation, the surface does not need to pass through \emph{any} of the vertices $v_0,\dots,v_i$, but the tangent plane at $v_0$ is perpendicular to $\hat n$. It is in principle possible to allow the fitted form to have a non-trivial tangent plane by including a term $bx+cy$ in \eqref{eq:quadratic:form}, however, we have not found this improves the quality of the results.

\subsubsection{On the non-robustness of vertex-based minimization}
The previous section might seem excessive effort for fitting a quadratic form. A more straightforward approach might be to simply minimize the distance of the form to the vertices:
\begin{equation*}
\tilde E[\bold{w}]=\sum_{j=1}^i  \prn{\begin{pmatrix}x_j\\y_j\end{pmatrix}^T
  \begin{pmatrix}
    d&e\\e&f
  \end{pmatrix}
 \begin{pmatrix}x_j\\y_j\end{pmatrix}
-z_j}^2.
\end{equation*}
It has been observed (for example by Meyer et al.~\cite{MeyerDesbrunSchroderBarr2002}) that this minimization can give poor results. Mathematically, minimizing this functional results in a linear system that may be under-resolved. For example, if $v_0,\dots,v_4$ comprise a rectangular pyramid, there are two trivial quadratic forms (and therefore an infinite number of them) that attain $\tilde E[\bold{w}]=0$ but have different mean curvatures. A vertex with three neighbors will always have a form that attains $\tilde E[\bold{w}]=0$, and the resulting quadratic forms often do not resemble the surface they are intended to approximate. This effect is drastic, and can result in a mean curvature of the opposite sign as the correct one, unstably causing a spike to emanate from the surface and grow indefinitely.

While we have no proof of robustness of our approach based on $E[\bold{w}]$, during the evolution of the ``Sculpture Flottante I'', with its 500,000 vertices to start, through merges and topology changes, until its demise not even once was linear system \eqref{eq:linear:system:w} poorly conditioned. By minimizing over the triangles (and not only the vertices) we also provide a robustness with regards to refining: if the triangulation is refined by adding a vertex on the edge which connects two neighbors of $v_0$, our functional $E[\bold{w}]$ remains unchanged while $\tilde E[\bold{w}]$ changes.

\subsection{Staleness}
\label{subsec:staleness}
When a triangulated surface evolves under mean curvature, the curvature at some vertices changes very quickly, while at other it changes slowly. Calculating the curvature at all the vertices at small time intervals can be computationally expensive, however, to get stable results  we must calculate it at short time intervals for vertices where curvature changes quickly.

To get around calculating the curvature for all the vertices at the rate set by the vertices who need the most frequent update, we recalculate the curvature for some vertices more frequently than for others. This is done by estimating for how long the curvature calculation remains relevant (``fresh'') at the same time that it is calculated. The ``expiration time'' we use for a vertex $v_0$ is given by
\begin{equation*}
\Delta t_{\text{exp}} = C_s \min_{j=1,\dots,n} \frac{\norm{v_j}}{\kappa}\;,
\end{equation*}
or, in words, some constant (we used $C_s = 0.5$) times the length of the minimal adjacent edge, divided by the velocity at which the the vertex moves, $\kappa$.

To avoid taking exceedingly small steps (stepping from one expired vertex to the next), once the curvature needs to be calculated at a vertex, we also recalculate all the vertices whose data is about to expire in the next time step. Then the largest step that can be taken (based on staleness) is performed. Of course, we also limit the step size by standard accuracy considerations, and make sure that we reach certain points in time to produce an output.

\subsection{Removal of Vertices}
\label{subsec:vertex_removal}
Since the size of the largest time step we are allowed to take is limited by the proximity of neighboring vertices, we take care to remove vertices that are too close to each other. We used three conditions as triggers to remove vertices:
\begin{itemize}
\item small edge, or
\item small angle, or
\item small face.
\end{itemize}
To remove a vertex, we merge two close vertices into a new one whose location is an average of the two. The ``freshness'' of the neighboring vertices (and the new vertex) is set to zero so that the curvature there be re-calculated immediately.

\begin{figure}[t]
\def\svgwidth{.98\textwidth}
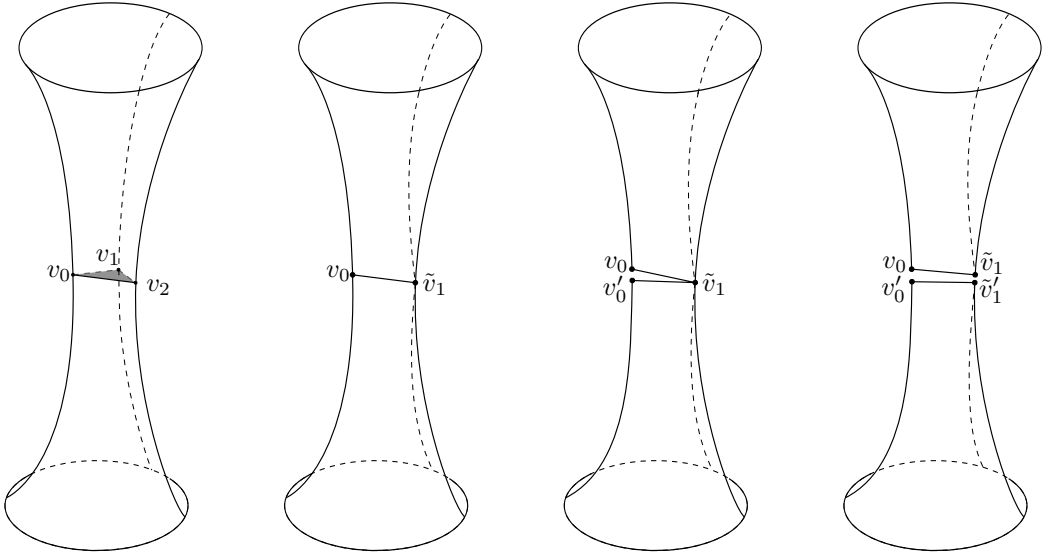
\caption{A topology change is initiated when two of the neighboring vertices $v_0$, $v_1$, and $v_2$, are to be merged into a vertex $\tilde v_1$, but the triangle $(v_0,\,v_1,\,v_2)$ (the shaded triangle), is not part of the surface. After the merge, the resulting surface is not homeomorphic to the 2-disk around $v_0$ and $\tilde v_1$. To resolve this, $v_0$ and then $\tilde v_1$ are duplicated (generating $v'_0$ and $\tilde v'_1$) and then connected to a subset of their original neighbors.}
\label{fig:topology:change}
\end{figure}

Special care needs to be taken to identify a required change in topology, and then make it happen. A change in topology is required when the neighborhood of a vertex is not homeomorphic to the 2-disk. As visualized in the stylized Fig.~\ref{fig:topology:change}, this can happen when two neighbors, $v_1$ and $v_2$, of a vertex, $v_0$, are merged, to $\tilde v_1$, but the surface does not contain the triangle $(v_0,\,v_1,\,v_2)$. This happens when the triangle forms the neck of a part that is about to pinch off. To identify this event, we examine the set of $v_0$'s neighbors, $\brk{v_1,\ldots,v_i}$ and attempt to order it so that for all $j$, $v_j$ and $v_{j+1}$ are neighbors. If this is not possible (because some vertex is needed in more then one place), there is a topology change. In this case, we find a subset of the vertices for which this is possible, and to each (perhaps not disjoint) set, we ``connect'' a unique copy of $v_0$. After splitting $v_0$, the same operation will be needed for $\tilde v_1$, as shown in Fig.~\ref{fig:topology:change}. These types of topology change treatments are similar to the ones presented by Mayer \cite{Mayer2001}, with the difference that in said paper surface diffusion flow is considered, and not mean curvature flow. It should also be remarked that unlike other types of geometric flows, for mean curvature flow the insertion of vertices is not required.

\section{The Characteristic Sounds of the Surface}
\label{sec:surface_sounds}
As described in Sect.~\ref{subsec:laplace_beltrami}, we use the surface wave equation \eqref{eq:wave_equation} to characterize the ``sound of the sculpture'' as it evolves under mean curvature flow. Since we are dealing with a surface that is evolving in time, $S(t)$, we also have a time-dependent family of Laplace-Beltrami operators. We denote this dependence on the surface by a subscript, i.e.~the function $u(t,x)$ which we use to characterize the ``sound of the sculpture'', as it evolves under mean curvature flow, satisfies
\begin{equation*}
u_{tt}(t,x)+F_s^2\LBO{S(t)} u(t,x)=0\quad \text{for all } t>0, \quad \text{ and } x\in S(t).
\end{equation*}
To simplify, we assume that the sound-generating vibrations have a much shorter time scale than that of the evolution of the surface itself, i.e.~$F_s^2\gg1$, and we allow ourselves to use a quasi-static approximation by which the surface itself is assumed to be constant between time steps. At every time step, we \emph{project} the solution onto the evolved surface and continue with the time evolution.

For a fixed surface, $S^n=S(t^n)$, solutions to
\begin{equation}
u_{tt}(t,x)+F_s^2\LBO{S^n} u(t,x)=0\;
\label{eq:wave:equation:t^n}
\end{equation}
can be found using separation of variables:
\begin{align}
u(t,x) &= e^{\pm i F_s\sqrt{\lambda}t} \phi(x)\;, \nonumber
\intertext{where}
\label{eq:eigenvalue:smooth}
-\LBO{S^n} \phi &= \lambda \phi\;.
\end{align}
Note that in the following, we consider the eigenvalues of the negative Laplace-Beltrami operator, since this is a positive semi-definite operator. The general solution to \eqref{eq:wave:equation:t^n} can be written using a linear combination of such simple terms:
\begin{equation}
u^n(t,x) = \sum_{k=1}^{\infty} \prn{c^n_{k,+}e^{ i F_s\sqrt{\lambda^n_k}(t-t^n)}+c^n_{k,-}e^{- i F_s\sqrt{\lambda^n_k}(t-t^n)}} \phi^n_k(x).
\label{eq:eignmode:sum}
\end{equation}
where $(\phi^n_k,\lambda^n_k)_{k=1}^\infty$ are all the eigenfunction/value pairs of $\LBO{S^n}$.

After evolving $u^n$ on the fixed surface $S^n$ for a time $\Delta t$, we project it and its time derivative, $u^n_t$, onto the surface $S^{n+1}$. The projected functions are used as initial condition for the evolution for another time step of $\Delta t$. From the resulting coefficients $c^n_{k,\pm}$ and the eigenvalues $\lambda^n_k$ we create the characteristic sound-wave of the evolving surface.

\subsection{Eigenfunctions and Eigenvalues}
To approximate the Laplace-Beltrami operator $\LBO{S^n}$ on the surface $S(t^n)$, we employ a ``cotan formula'' \eqref{eq:discrete_Laplace_Beltrami} with the specific weights \eqref{eq:LB_weights}, as it is given in \cite{BobenkoSpringborn2007}. This discretization replaces the infinite dimensional eigenvalue problem \eqref{eq:eigenvalue:smooth} by a matrix eigenvalue problem
\begin{equation*}
A^n \varphi = \lambda \varphi\;,
\end{equation*}
where $\varphi = (\varphi(v_1),\dots,\varphi(v_{|S^n|}))^T$, and the entries of the matrix $A^n$ are given by the weights in the cotan formula, as
\begin{equation}
a^n_{ii} = \sum_{j} w_{ij}\;\forall\,i
\quad\text{and}\quad
a^n_{ij} = -w_{ij}\;\forall\,i\neq j\;,
\label{eq:LBO:discretization}
\end{equation}
which, of course, depend on the specific triangulation of $S^n$ at time $t^n$.
Consequently, the matrix $A^n$ is a large (initially about $500,000\times500,000$), sparse, symmetric matrix, and we need to find its smallest $N=50$ eigenvalues and corresponding eigenvectors. We achieve this task by an implicitly restarted Arnoldi method \cite{LehoucqSorensen1996} implemented in Matlab's \texttt{eigs}. The symmetry of $A$ arises from the self-adjointness of $\LBO{}$ and from the symmetry of the cotan formula. It implies that the resulting eigenfunctions are orthogonal (under summation over the vertices). They are normalized so that
\begin{equation*}
\varphi^n_k\cdot\varphi^n_l=\delta_{kl}\;.
\end{equation*}
In this notation, $\varphi^n_k$ denotes the $k$-th eigenfunction of the discrete Laplace-Beltrami operator, discretized according to \eqref{eq:LBO:discretization} on the triangulation at time $t=t^n$.

\begin{figure}
\begin{minipage}[b]{.48\textwidth}
\includegraphics[width=\textwidth]{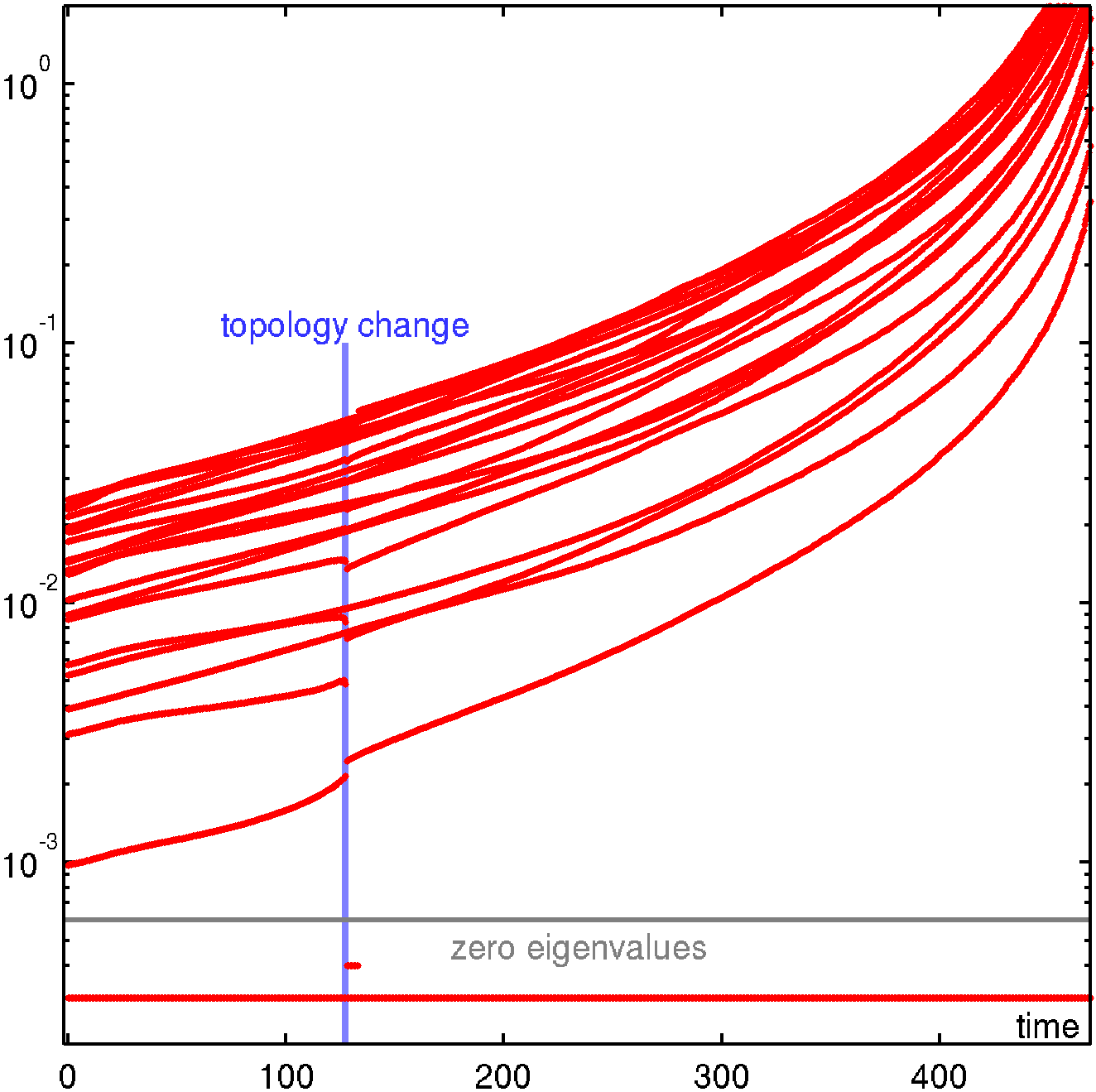}
\caption{Temporal evolution of the eigenvalues of the sculpture's top piece.}
\label{fig:eigenvalues}
\end{minipage}
\hfill
\begin{minipage}[b]{.48\textwidth}
\includegraphics[width=\textwidth]{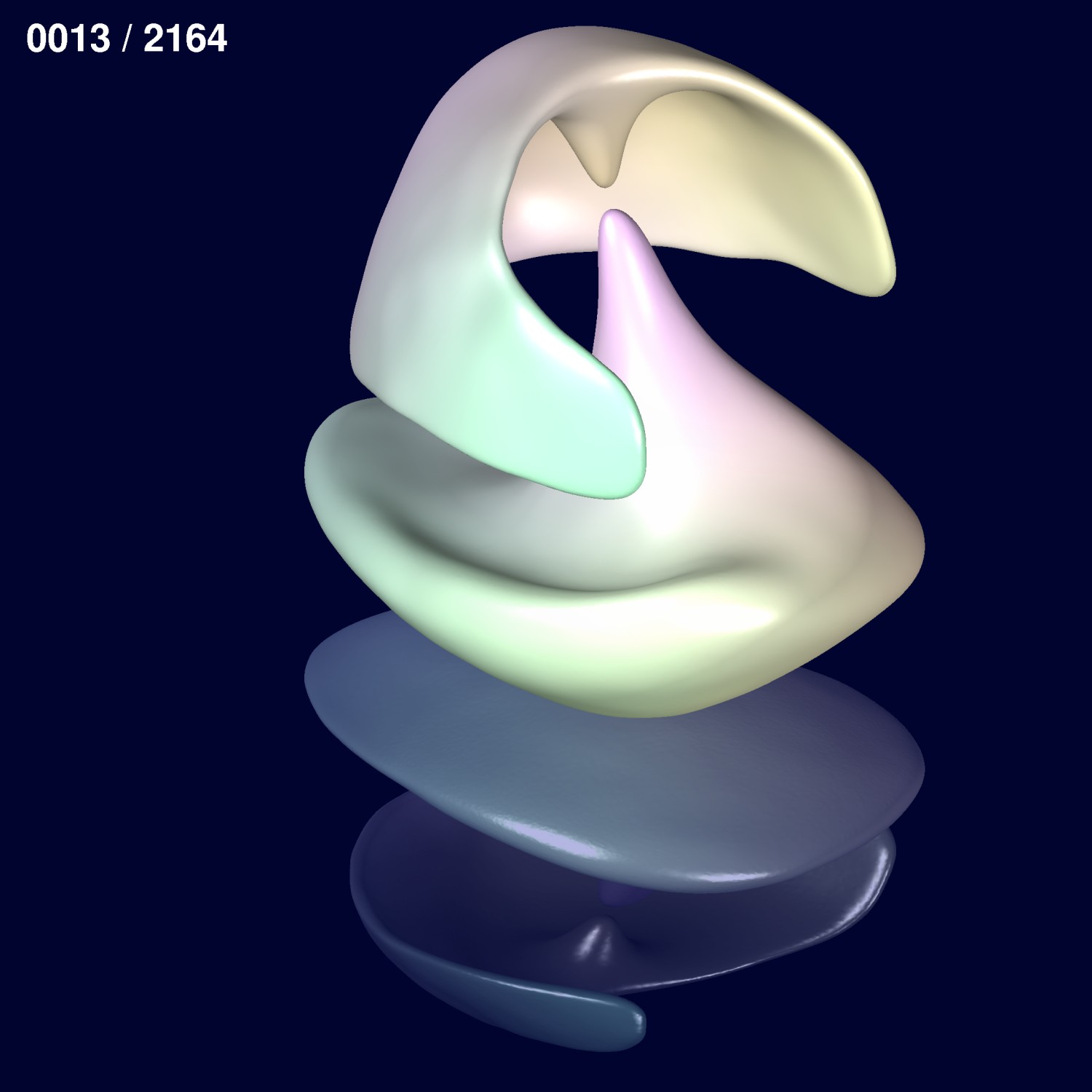}
\caption{A frame of the final video of the ``Evolving Floating Sculpture''.}
\label{fig:video_frame}
\end{minipage}
\end{figure}

The temporal evolution of the $N$ lowest eigenvalues of the sculpture's top part, as it evolves under mean curvature flow, is shown in Fig.~\ref{fig:eigenvalues}. One can see how the frequencies increase as the surface shrinks. Moreover, one can observe the discontinuity in the eigenvalues due to the pinch-off event (shown in Fig.~\ref{fig:animation}, between panels (c) and (d)). The eigenvalues that equal zero correspond directly to the connected components of the surface. The brief appearance of a second connected component after the pinch-off is echoed by a short segment in the ``zero eigenvalues'' zone.

\subsection{Initial Condition and the Evolution During a Time Step}
Equation~\eqref{eq:wave:equation:t^n} requires initial conditions for both $u$ and $u_t$. For all $n>0$ we use the projection from the end of the previous step as the initial conditions for the current step. The initial condition at $n=0$ is chosen to attain 1 on a single vertex (chosen manually), and to vanish on the rest. The time derivative, $u_t$, is chosen to vanish identically at $n=0$ so that we can visualize this initial condition as ``plucking'' the mesh at one vertex.

At the beginning of each step the values of $c^n_{k,\pm}$ are found from the values of $u^n$ and $u^n_t$ by solving a sequence of small linear systems:
\begin{equation*}
\begin{pmatrix}
\varphi^n_m\cdot u^n(t,\cdot)\\
\varphi^n_m\cdot u^n_t(t,\cdot)
\end{pmatrix}=
\begin{pmatrix}
1&1\\
iF_s\sqrt{\lambda^n_m}  &-iF_s\sqrt{\lambda^n_m}
\end{pmatrix}
\begin{pmatrix}
c^n_{m,+}(t)\\
c^n_{m,-}(t)
\end{pmatrix}\;.
\end{equation*}
This system is found from \eqref{eq:eignmode:sum} by differentiating with respect to $t$ and then taking the inner product with the eigenfunction $\varphi^n_m$.

Once the coefficients of the solution at time $t=t^n$ have been found, the solution and its time derivative are evolved to $t=t^{n+1}$ using the general solution \eqref{eq:eignmode:sum}:
\begin{align*}
\tilde u^n(t^{n+1},v)&=\sum_{k=1}^{N} \prn{c^n_{k,+}e^{ i F_s\sqrt{\lambda^n_k}\Delta t}+c^n_{k,-}e^{- i F_s\sqrt{\lambda^n_k}\Delta t}} \varphi^n_k(v)\;, \\
\tilde u_t^n(t^{n+1},v)&=\sum_{k=1}^{N} iF_s\prn{c^n_{k,+}e^{ i F_s\sqrt{\lambda^n_k}\Delta t}-c^n_{k,-}e^{- i F_s\sqrt{\lambda^n_k}\Delta t}} \varphi^n_k(v)\;,
\end{align*}
where we write $\tilde u$ and not $u$ to signify that it was not yet projected on $S(t^{n+1})$.

\subsection{Projection}
Before we can take another time step we must project the function and its time derivative from the surface $S(t^n)$ onto $S(t^{n+1})$.
During the evolution, we keep track of the merging of vertices, so that for every vertex $v^{n+1}_k$ at time $t^{n+1}$ we can identify a ``parent vertex'' $v^n_l$ at time $t^n$ with the property that the evolution of the surface takes $v^n_l$ onto  $v^{n+1}_k$ by evolution and merges. This ``parent image'' is used for the projection: the value of a function $f(v^{n+1}_k)$ is taken to be $\tilde f(v^n_l)$.

\subsection{The Wave-Form}
Once we have found the coefficients $c^n_{k,\pm}$ for the evolution of the sculpture until it vanishes, we construct the characteristic sounds of the evolution. We approximate continuous functions $\Lambda_k(t)$ and $C_k(t)$ by a linear interpolation of the $\lambda^n_k$ and the $c^n_{k,+}$. Using these, we generate the waveform $W(t)$ as the sum
\begin{equation*}
W(t) = \sum_{k=1}^N C^n_k(t) \sin\!\prn{F_s \int_0^t \Lambda_k(t')\ud{t'} }\;,
\end{equation*}
which is then sampled at a high resolution to create a \texttt{.wav} file.

\section{The Transformation of the Computational Results into Pieces of Art}
\label{sec:results_to_art}
Having created the computational results from the mean curvature flow of the sculpture's surface (see Sect.~\ref{sec:geometry_evolution}), and having an approach to compute its characteristic sounds (see Sect.~\ref{sec:surface_sounds}), we now describe the crucial steps of assembling these data into pieces of art.

\begin{figure}
\begin{minipage}[b]{.24\textwidth}
\includegraphics[width=\textwidth]{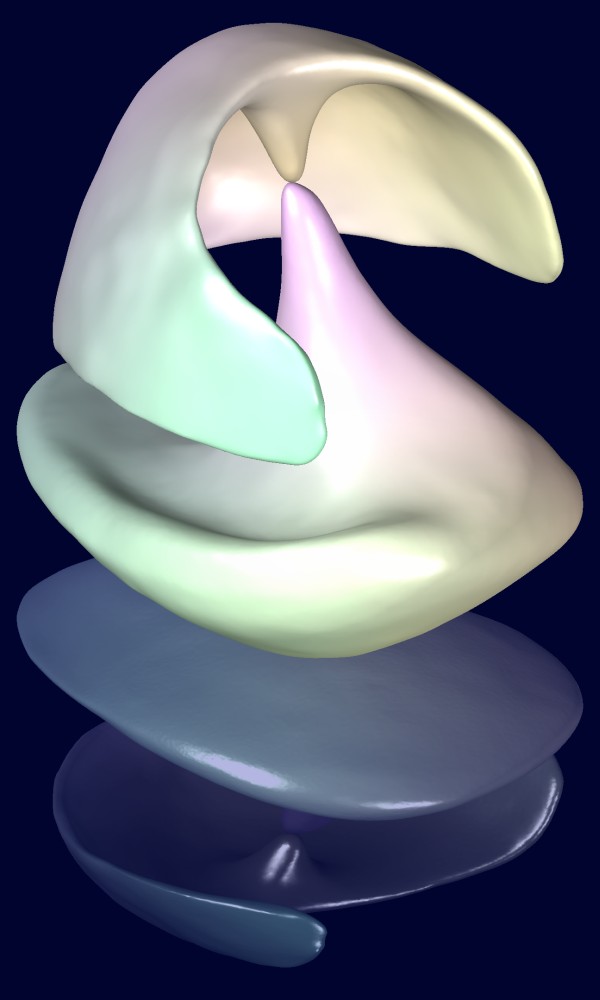}
\hspace{-9.5em}\parbox{9em}{\flushright{\textcolor{white}{initial surface (a)}}\\~}
\end{minipage}
\hfill
\begin{minipage}[b]{.24\textwidth}
\includegraphics[width=\textwidth]{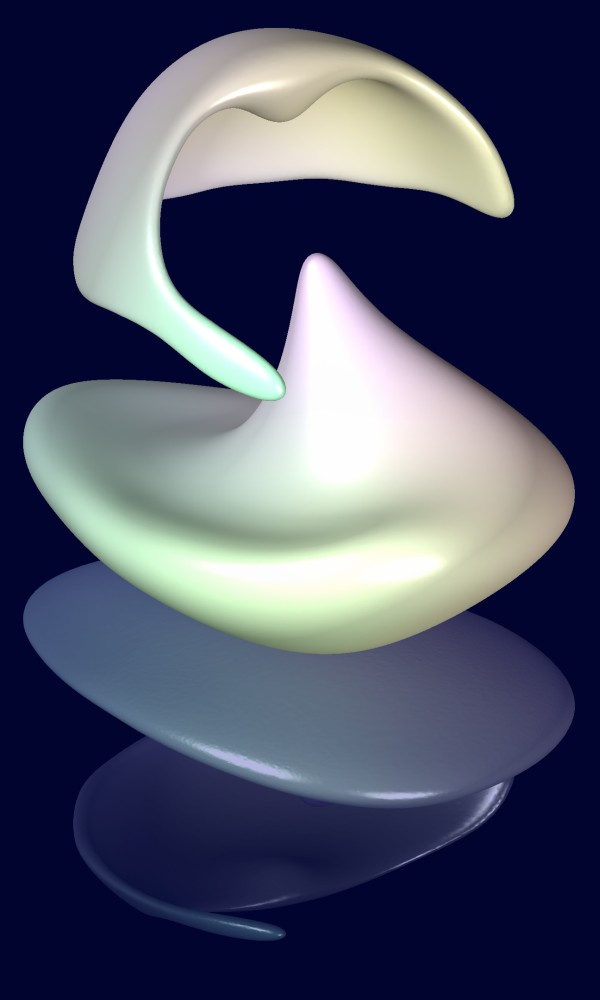}
\hspace{-9.5em}\parbox{9em}{\flushright{\textcolor{white}{(b)}}\\~}
\end{minipage}
\hfill
\begin{minipage}[b]{.24\textwidth}
\includegraphics[width=\textwidth]{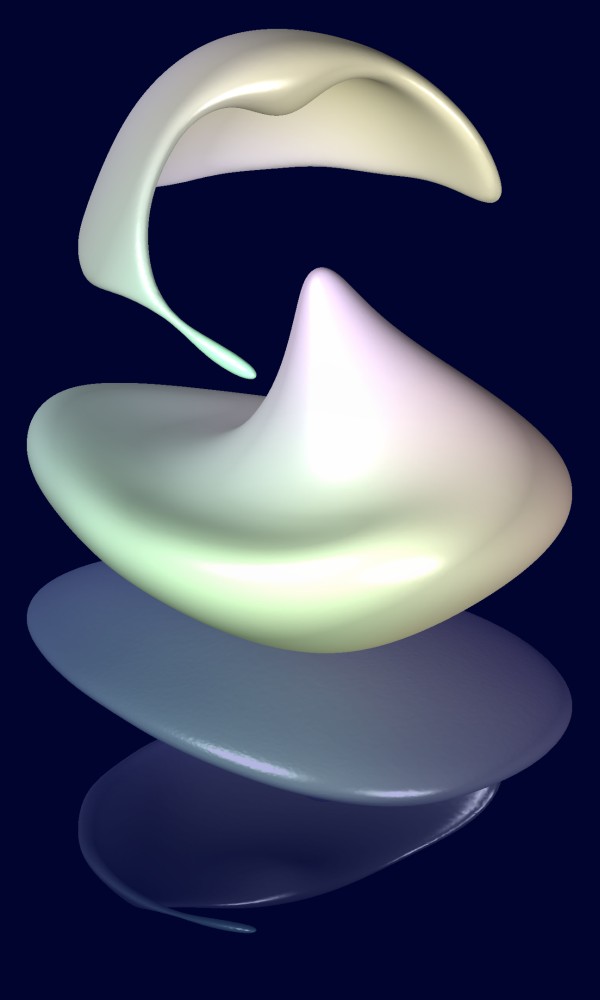}
\hspace{-9.5em}\parbox{9em}{\flushright{\textcolor{white}{before pinch-off (c)}}\\~}
\end{minipage}
\hfill
\begin{minipage}[b]{.24\textwidth}
\includegraphics[width=\textwidth]{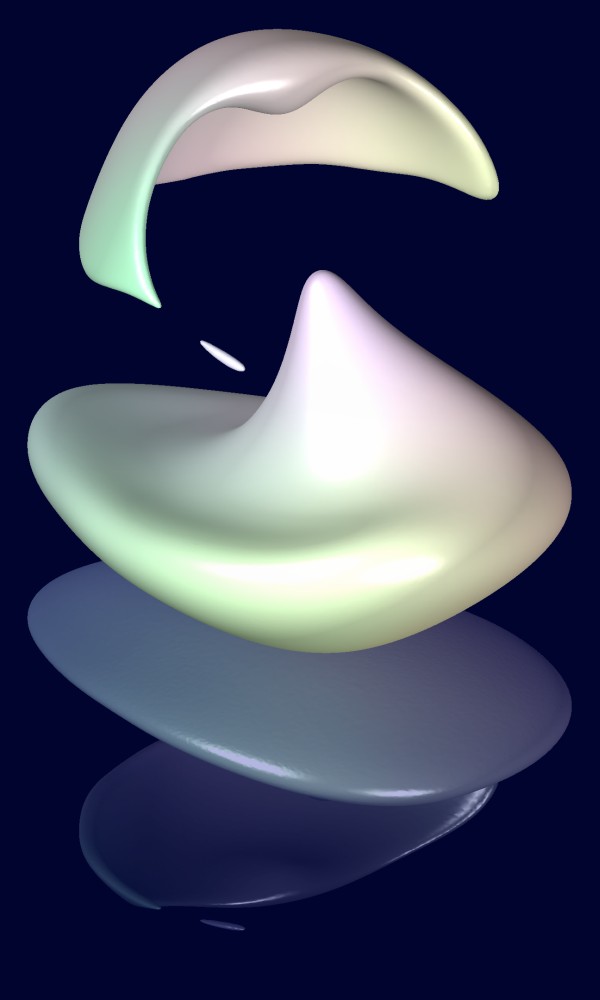}
\hspace{-9.5em}\parbox{9em}{\flushright{\textcolor{white}{after pinch-off (d)}}\\~}
\end{minipage}

\vspace{-.3em}
\begin{minipage}[b]{.24\textwidth}
\includegraphics[width=\textwidth]{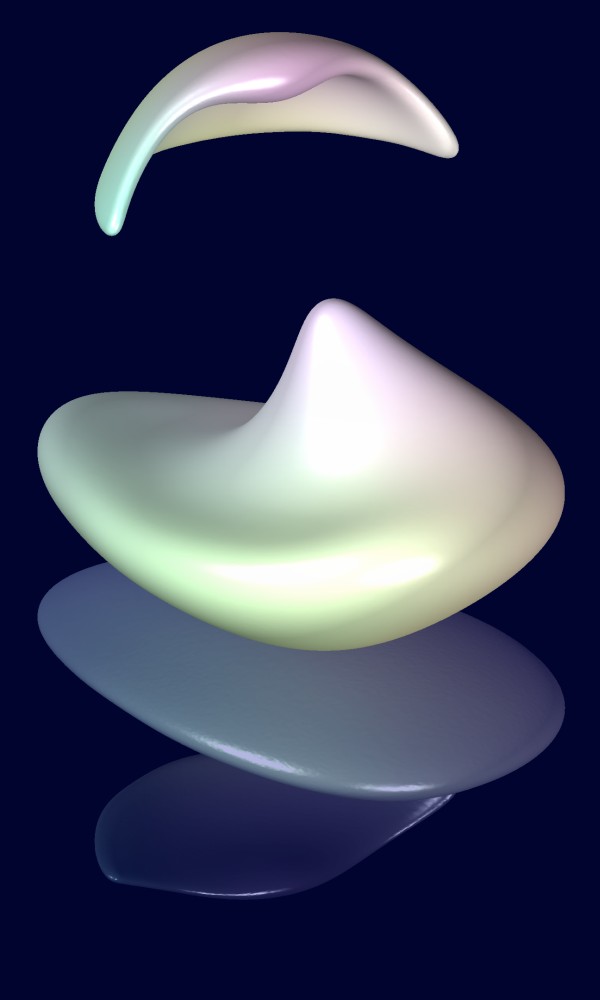}
\hspace{-9.5em}\parbox{9em}{\flushright{\textcolor{white}{(e)}}\\~}
\end{minipage}
\hfill
\begin{minipage}[b]{.24\textwidth}
\includegraphics[width=\textwidth]{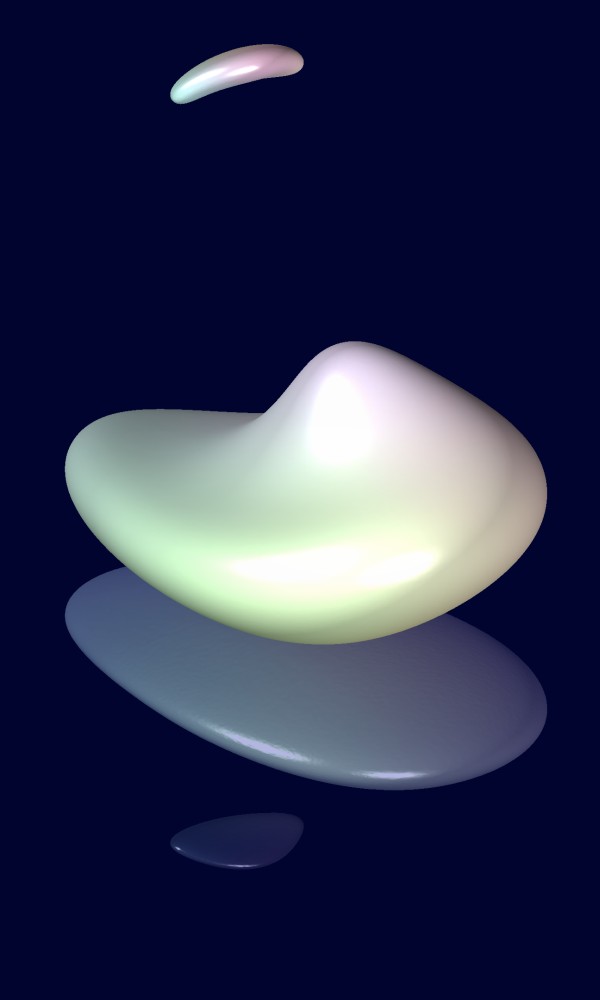}
\hspace{-9.5em}\parbox{9em}{\flushright{\textcolor{white}{top almost gone (f)}}\\~}
\end{minipage}
\hfill
\begin{minipage}[b]{.24\textwidth}
\includegraphics[width=\textwidth]{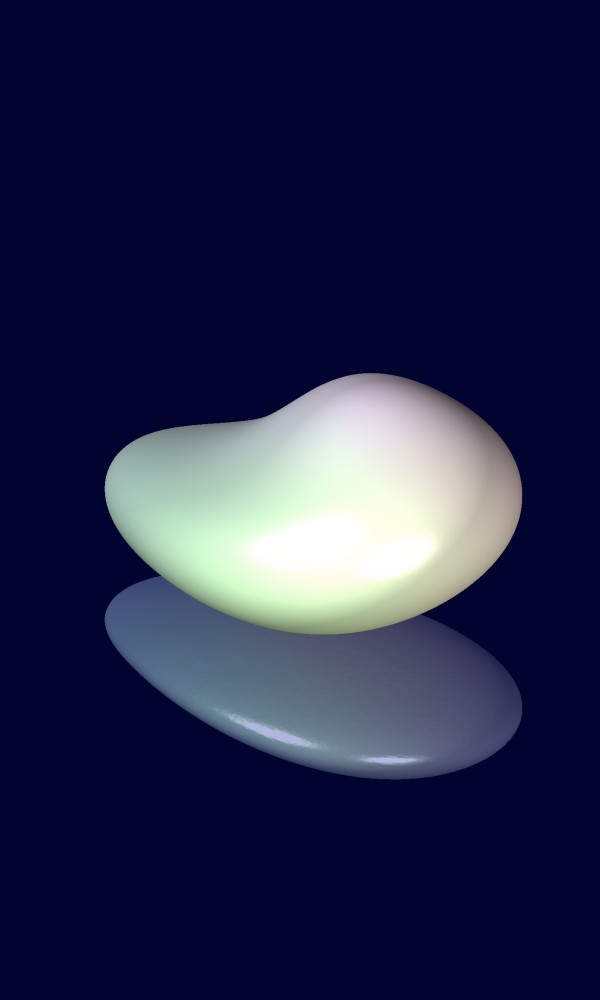}
\hspace{-9.5em}\parbox{9em}{\flushright{\textcolor{white}{top has vanished (g)}}\\~}
\end{minipage}
\hfill
\begin{minipage}[b]{.24\textwidth}
\includegraphics[width=\textwidth]{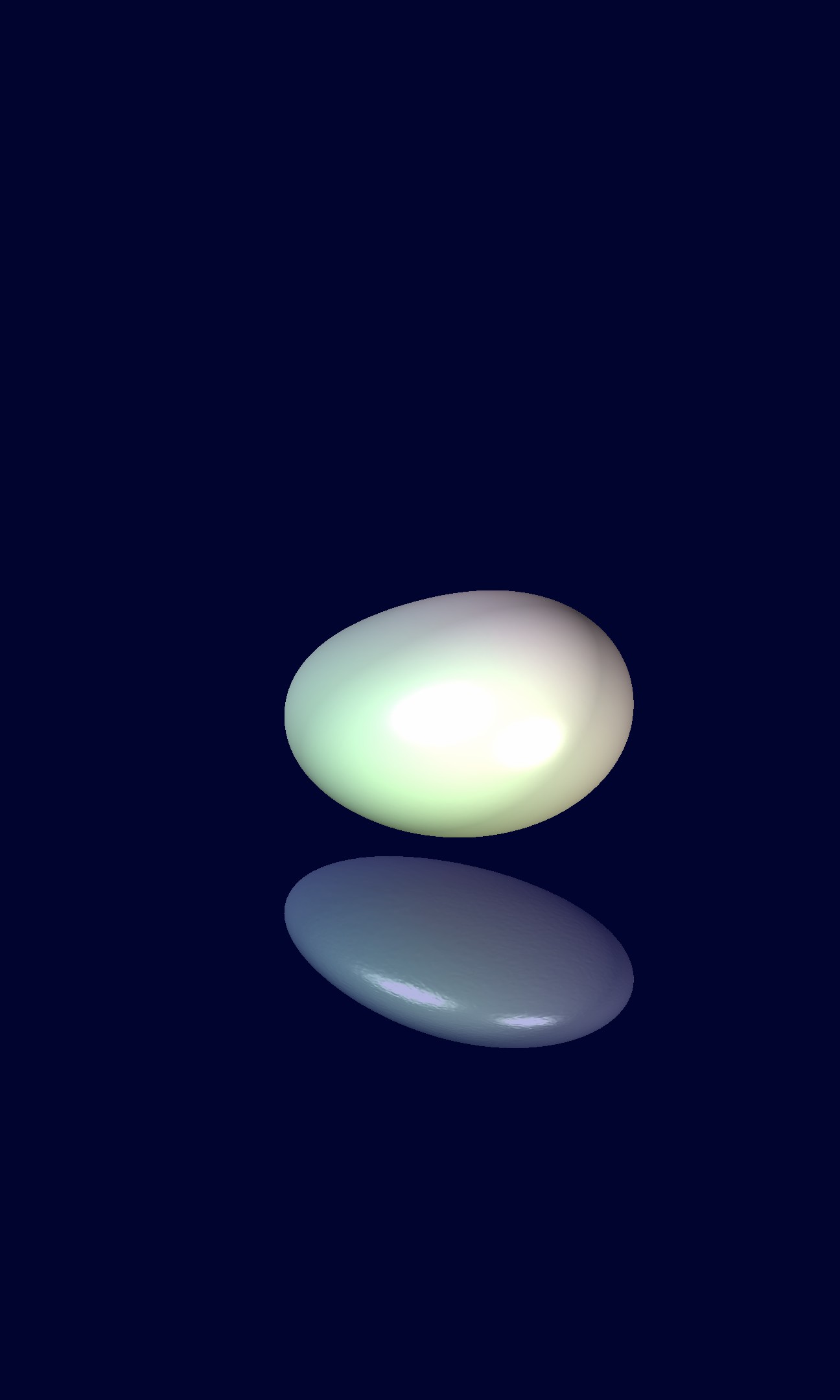}
\hspace{-9.5em}\parbox{9em}{\flushright{\textcolor{white}{(h)}}\\~}
\end{minipage}
\caption{The evolution of the sculpture's surface under mean curvature flow. From top left to bottom right: (a) the initial configuration, and the surface after (b)~100~, (c)~125~, (d)~130~, (e)~200~, (f)~400~, (g)~700~, (h)~1500 time units. Between (c) and (d) a small piece is pinched-off the top part.}
\label{fig:animation}
\end{figure}

\subsection{The Video}
A video is created that shows the evolution of the sculpture's surface under mean curvature flow, from its initial state, until it vanishes completely. The triangulation is rendered using a Phong shading illumination model \cite{Phong1975}, involving ambient, diffuse, and specular reflections of two light sources. At any instance in time, each connected component receives a subtle coloring based on its three lowest (non-constant) eigenfunctions, i.e.~the ones that correspond to the three smallest, positive eigenvalues (there is a zero eigenvalue for every connected component, corresponding to a constant eigenfunction on that component). Specifically, the lowest, second-lowest, and third-lowest frequency eigenfunctions affect the red, green, and blue of the surface's RGB channels, respectively. As described in \cite{Levy2006, ReuterBiasottiGiorgiPataneSpagnuolo2009}, the resulting coloring of the surface's connected components adds some characteristic geometric information to the visualization.

In discussion with the artist Jane Philbrick, a number of artistic features are added. First, the evolving surface is placed on top of a blue plane, and a mirrored, blue-shaded copy of the surface is placed below the plane, thus imitating reflection in water, as it would occur with the original ``Sculpture Flottante I'' (see Fig.~\ref{fig:sculpture_flottante}). This water-like effect is further amplified by adding small ripples on the reflection. Second, an animation sequence is added to precede the actual mean curvature flow, in which the digital version of the ``Sculpture Flottante I'' is introduced to the viewer. In this sequence, the sculpture is initially approached by a moving camera, and then rotated to allow a $360^\circ$ view of its shape. Third, a counter is added that clocks the pseudo-time that has passed from the onset of the mean curvature flow, providing a reference for the various time scales that the flow possesses: the fast contraction of high-curvature regions on the one hand, and the slow shrinking of almost flat areas on the other hand. A frame of the video is shown in Fig.~\ref{fig:video_frame}. At the instance shown, the mean curvature flow has just started, and created a small separation between the top and the bottom pieces. A sequence of eight representative snapshots is given in Fig.~\ref{fig:animation}. The initial configuration is shown in panel (a). Between panels (b) and (c), a small pinch-off event happens (see Fig.~\ref{fig:eigenvalues} for the corresponding behavior of the eigenvalues). Between panels (f) and (g), the top piece vanishes, and finally in panel (h) the base piece is approaching the last stages of its existence. The full video can be downloaded on the project website \cite{FloatingSculptureProject}.

\subsection{The Sculpture Song}
The approach described in Sect.~\ref{sec:surface_sounds} yields the temporal evolution of an initial ``plucking'' of a single point in one connected component of the surface, as it evolves under mean curvature flow. When, where, and how frequently to excite the pieces of the surface to create an interesting sound piece is an artistic decision. As such, we created a number of sound files, arising from different initial conditions and provided them to the artist, who then mixed these pieces together to create the final ``Sculpture Song''. The resulting sound file is available for downloaded on the project website \cite{FloatingSculptureProject}.

\section*{Acknowledgments}
The authors would like to express their gratitude to Jane Philbrick (MIT Center for Advanced Visual Studies), who approached us mathematicians with an open mind for a joint project, for many inspiring discussions, and for providing the stage for our computational projects to turn into pieces of art. We also would like to thank Enno Lenzmann (now University of Copenhagen), who established the first contact between us and Jane Philbrick, for many helpful comments during the creation of the visualization. The authors would also like to acknowledge the support by the National Science Foundation. B. Seibold was supported through grant DMS--0813648, and Y. Farjoun was supported through grant DMS--0703937. Moreover, Y. Farjoun wishes to acknowledge partial support by the Spanish Ministry of Science and Innovation through grant FIS2008-04921-C02-01, and B. Seibold wishes to acknowledge partial support by NSF through grants DMS--1007899 and DMS--1115278.

\bibliographystyle{plain}
\bibliography{references_complete}

\end{document}

%% file: TopologyChange.pdf_tex
\begingroup%
  \makeatletter%
  \providecommand\color[2][]{%
    \errmessage{(Inkscape) Color is used for the text in Inkscape, but the package 'color.sty' is not loaded}%
    \renewcommand\color[2][]{}%
  }%
  \providecommand\transparent[1]{%
    \errmessage{(Inkscape) Transparency is used (non-zero) for the text in Inkscape, but the package 'transparent.sty' is not loaded}%
    \renewcommand\transparent[1]{}%
  }%
  \providecommand\rotatebox[2]{#2}%
  \ifx\svgwidth\undefined%
    \setlength{\unitlength}{992.90166016bp}%
    \ifx\svgscale\undefined%
      \relax%
    \else%
      \setlength{\unitlength}{\unitlength * \real{\svgscale}}%
    \fi%
  \else%
    \setlength{\unitlength}{\svgwidth}%
  \fi%
  \global\let\svgwidth\undefined%
  \global\let\svgscale\undefined%
  \makeatother%
  \begin{picture}(1,0.55640212)%
    \put(0,0){\includegraphics[width=\unitlength]{TopologyChange.pdf}}%
    \put(0.06699358,0.27771845){\color[rgb]{0,0,0}\makebox(0,0)[lb]{\smash{$v_0$}}}%
    \put(0.15774709,0.26511108){\color[rgb]{0,0,0}\makebox(0,0)[lb]{\smash{$v_2$}}}%
    \put(0.11174657,0.29216228){\color[rgb]{0,0,0}\makebox(0,0)[lb]{\smash{$v_1$}}}%
    \put(0.32121364,0.27771845){\color[rgb]{0,0,0}\makebox(0,0)[lb]{\smash{$v_0$}}}%
    \put(0.40874416,0.26511108){\color[rgb]{0,0,0}\makebox(0,0)[lb]{\smash{$\tilde v_1$}}}%
    \put(0.57272988,0.26140937){\color[rgb]{0,0,0}\makebox(0,0)[lb]{\smash{$v'_0$}}}%
    \put(0.66296442,0.26511103){\color[rgb]{0,0,0}\makebox(0,0)[lb]{\smash{$\tilde v_1$}}}%
    \put(0.57366788,0.28619584){\color[rgb]{0,0,0}\makebox(0,0)[lb]{\smash{$v_0$}}}%
    \put(0.82589026,0.26140937){\color[rgb]{0,0,0}\makebox(0,0)[lb]{\smash{$v'_0$}}}%
    \put(0.91545259,0.2844608){\color[rgb]{0,0,0}\makebox(0,0)[lb]{\smash{$\tilde v_1$}}}%
    \put(0.82685601,0.28619583){\color[rgb]{0,0,0}\makebox(0,0)[lb]{\smash{$v_0$}}}%
    \put(0.91469488,0.2599205){\color[rgb]{0,0,0}\makebox(0,0)[lb]{\smash{$\tilde v'_1$}}}%
  \end{picture}%
\endgroup%